\documentclass{amsart}

\usepackage{graphicx}
\usepackage{amsmath}
\usepackage{amsfonts}
\usepackage{amssymb}

\newtheorem{theorem}{Theorem}[section]

\newtheorem{proposition}[theorem]{Proposition}

\theoremstyle{definition}

\theoremstyle{remark}

\newcommand{\be}{\begin{equation}}
\newcommand{\ee}{\end{equation}}
\newcommand{\ben}{\begin{equation*}}
\newcommand{\een}{\end{equation*}}

\newcommand{\bma}{\begin{pmatrix}}
\newcommand{\ema}{\end{pmatrix}}
\newcommand{\ba}{\begin{array}}
\newcommand{\ea}{\end{array}}
\newcommand{\beq}{\begin{eqnarray}}
\newcommand{\eeq}{\end{eqnarray}}
\newcommand{\beqn}{\begin{eqnarray*}}
\newcommand{\eeqn}{\end{eqnarray*}}
\newcommand{\balign}{\begin{align}}
\newcommand{\ealign}{\end{align}}
\newcommand{\bal}{\begin{aligned}}
\newcommand{\eal}{\end{aligned}}

\newcommand{\der}{{\rm d}}
\newcommand{\wt}[1]{\widetilde{#1}}
\newcommand{\gt}[1]{\mathfrak{#1}}
\newcommand{\w}{\wedge}

\begin{document}
\title{Third order ODEs and four-dimensional split signature Einstein metrics} 
\vskip 1.truecm \author{Micha\l~ Godli\'nski} \address{Instytut Fizyki Teoretycznej,
Uniwersytet Warszawski, ul. Ho\.za 69, Warszawa, Poland}
\email{godlinsk@fuw.edu.pl} \thanks{This research was supported by
the KBN grant 2 P03B 12724}

\author{Pawe\l~ Nurowski} \address{Instytut Fizyki Teoretycznej,
Uniwersytet Warszawski, ul. Hoza 69, Warszawa, Poland}
\email{nurowski@fuw.edu.pl} \thanks{During the preparation of this
article P. N. was a memeber of the VW Junior Research Group ``Special
Geometries in Mathematical Physiscs'' at Humboldt University in
Berlin.}

\date{\today}
\begin{abstract}
We construct a family of split signature Einstein metrics in four dimensions, 
corresponding to particular classes of third order ODEs considered 
modulo fiber preserving transformations of variables.
\end{abstract}
\maketitle

\allowdisplaybreaks

\section{Introduction}
Our starting point is a 3rd order ordinary differential equation (ODE)
\begin{equation}
y'''=F(x,y,y',y'')\label{r1},
\end{equation}
for a real function $y=y(x)$. Here $F=F(x,y,p,q)$ is 
a sufficiently smooth real function of four real variables
$(x,y,p=y',q=y'')$.

Given another 3rd order ODE 
\be
\bar{y}'''=\bar{F}(\bar{x},\bar{y},\bar{y}',\bar{y}'')\label{r1p}
\ee
it is
often convenient to know whether there exists a suitable transformation of
variables $(x,y,p,q)\to(\bar{x},\bar{y},\bar{p},\bar{q})$ which brings 
(\ref{r1p}) to (\ref{r1}). 
Several types of such transformations are of particular importance.
Here we consider fiber preserving (f.p.) transformations, which are of the form 
\be
\bar{x}=\bar{x}(x),\quad\quad\quad\quad\quad\quad\quad\quad\quad\bar{y}=\bar{y}(x,y).
 \label{r2}
\ee
We say that two 3rd order ODEs, (\ref{r1}) and \eqref{r1p}, are
(locally) f.p. equivalent iff there exists a (local) f.p. transformation \eqref{r2},
which brings (\ref{r1p}) to (\ref{r1}). The task of finding neccessary
and sufficient conditions for ODEs (\ref{r1}) and (\ref{r1p}) to be 
(locally) f.p. equivalent, is called a f.p. equivalence problem for 3rd order
ODEs. In the cases of (more general) point transformations and contact transformations, this problem was
studied and solved by  E. Cartan \cite{Cartanode} and S-S. Chern \cite{Chern} in the years 1939-1941. 
The interest in these studies has been recently revived due to the
fact that 
important 
equivalence classes of 3rd order ODEs naturally define 3-dimensional conformal 
Lorentzian structures including Einstein-Weyl structures. This makes these equivalence
problems aplicable not only to differential geometry but also to
the theory of integrable systems and General Relativity \cite{Newman,Nurconf,Tod}.

In this paper we show how to construct 4-dimensional split signature Einstein metrics,
starting from particular ODEs of 3rd order. 
We formulate the problem of f.p. equivalence in terms of differential forms. 
Invoking Cartan's equivalence method, we construct a 6-d manifold with a distinguished coframe on it,
which encodes all information about original equivalence problem. 
For specific types of the ODEs, the class of Einstein metrics can be explicitly constructed from this coframe. 
This result is a byproduct of the full solution 
of the f.p. equivalence problem, that will be described in \cite{Nurgodl}.

We acknowledge that all our calculations were checked by the
independent use of the two symbolic calculations programs: Maple and 
Mathematica.
\section{Third order ODE and Cartan's method}
Following Cartan and Chern, we rewrite \eqref{r1}, using 1-forms 
\be\bal
 \omega^1&=\der y -p\der x, \\
 \omega^2&=\der p- q\der x, \\
 \omega^3&=\der q- F(x,y,p,q)\der x,\\     
 \omega^4&=\der x. \\
\eal\label{omega}
\ee
These are defined on the second jet space $\mathcal{J}^2$ locally parametrized by $(x,\,y,\,p,\,q)$. 
Each solution $y=f(x)$ of \eqref{r1} is fully described by the two conditions: forms $\omega^1$, $\omega^2$, $\omega^3$ vanish 
on a curve $(t,f(t),f'(t),f''(t))$ and, as this defines a solution up to transformations of $x$, 
$\omega^4=\der t$ on this curve. Suppose now, that  equation \eqref{r1} undergoes fiber preserving transformations \eqref{r2}. 
Then the forms \eqref{omega} transform by
\be\bal
 \omega^1&\to \bar{\omega}^1=\alpha\omega^1, \\
 \omega^2&\to \bar{\omega}^2=\beta(\omega^2+\gamma\omega^1),  \\
 \omega^3&\to \bar{\omega}^3=\epsilon(\omega^3+\eta\omega^2+\varkappa\omega^1),\\     
 \omega^4&\to \bar{\omega}^4=\lambda\omega^4,  \\
\eal\label{tranom}
\ee
where functions $\alpha,\,\beta,\,\gamma,\,\epsilon,\,\eta,\,\varkappa,\,\lambda$ are defined on 
$\mathcal{J}^2$, satisfy $\alpha\beta\epsilon\lambda\neq 0$ and are determined by a particular
choice of transformation \eqref{r2}. 
A fiber preserving equivalence class of ODEs is described by forms 
\eqref{omega} defined up to transformations \eqref{tranom}. Equations \eqref{r1} and \eqref{r1p} 
are f.p. equivalent, iff their corresponding forms $(\omega^i)$ and $(\bar{\omega}^j)$ are related as above.

We now apply Cartan's equivalence method 
\cite{Olver,Sternberg}. Its key idea is to enlarge the space $\mathcal{J}^2$ 
to a new manifold $\mathcal{\wt{P}}$, on which functions $\alpha,\,\beta,\,\gamma,\,\epsilon,\,\eta,\,\varkappa,\,\lambda$ 
are additional coordinates. The coframe $(\omega^i)$ defined up to transformations \eqref{tranom}, 
is now replaced by a set of four well defined 1-forms 
\ben\bal
 \theta^1&=\alpha\omega^1, \\
 \theta^2&=\beta(\omega^2+\gamma\omega^1),  \\
 \theta^3&=\epsilon(\omega^3+\eta\omega^2+\varkappa\omega^1),\\     
 \theta^4&=\lambda\omega^4  \\
\eal
\een
on $\mathcal{\wt{P}}$. If, in addition, the following f.p. invariant
condition \cite{Godlinski,Grebot}
$$F_{qq}\neq 0$$
is satisfied then, there is a geometrically distinguished way of 
choosing five parameters
$\beta,\,\epsilon,\,\eta,\,\varkappa,\,\lambda$ to be functions of
$(x,y,p,q,\alpha,\gamma)$. Then, on a 6-dimensional manifold
$\mathcal{P}$ parametrized by $(x,y,p,q,\alpha,\gamma)$ Cartan's
method give a way of supplementing the well defined four 1-forms $(\theta^i)$
with two other 1-forms $\Omega^1$, $\Omega^2$ so that the set
$(\theta^1,\theta^2,\theta^3,\theta^4,\Omega^1,\Omega^2)$ constitutes
a rigid coframe on $\mathcal{P}$. According to the theory of G-structures \cite{Kobayashi,Sternberg}, 
all information about a f.p. equivalence class of equation \eqref{r1} satisfying $F_{qq}\neq 0$ is encoded in the coframe 
$(\theta^1$,$\theta^2$,$\theta^3$,$\theta^4$,$\Omega^1$,$\Omega^2)$. 
Two equations \eqref{r1} and \eqref{r1p} are f.p. equivalent, iff there exists a diffeomorphism 
$\psi:\mathcal{P}\to\mathcal{\bar{P}}$, such that
$\psi^{\ast}\bar{\theta}^i=\theta^i$,
$\psi^{\ast}\bar{\Omega}^A=\Omega^A$, where 
$i=1,2,3,4$ and $A=1,2$.
The procedure of constructing manifold $\mathcal{P}$ and the coframe $(\theta^i,\Omega^A)$ is explained in details in 
\cite{Olver,Sternberg} for a general case and in \cite{Godlinski,Nurgodl} for this specific problem. Here we omit the details of this procedure, 
summarizing the results on f.p. equivalence problem in the following theorem.
\begin{theorem}\label{tw2}
A 3rd order ODE $y'''=F(x,y,y',y'')$,  satisfying $F_{qq}\neq 0$, 
considered modulo fiber preserving transformations of variables, uniquely defines 
a 6-dimensional manifold $\mathcal{P}$, and an invariant coframe  $(\theta^1,\theta^2,\theta^3,\theta^4,\Omega^1,\Omega^2)$ on it. 
In local coordinates $(x,y,p=y',q=y'',\alpha,\gamma)$ this coframe is given by
 \beq\label{theta}
  &\theta^1=\alpha\omega^1,\nonumber \\
  &\theta^2=\tfrac{1}{6}F_{qq}(\omega^2+\gamma\omega^1),\nonumber \\
  &\theta^3=\tfrac{1}{36\alpha}F_{qq}(\omega^3+(\gamma-\tfrac{1}{3})F_{q}\omega^2+(\tfrac{1}{2}\gamma^2+K)\omega^1),
   \nonumber\\
  &\theta^4=\tfrac{6\alpha}{F_{qq}}\omega^4,\nonumber \\
  &\Omega^1=\tfrac{1}{F_{qq}}(-F_{qqq}\gamma^2+(\tfrac{2}{3}F_{qqq}F_{q}+\tfrac{1}{3}F_{qq}^2+2F_{qqp})\gamma\nonumber\\
     &+F_{qq}K_{q}+2F_{qqq}K-2F_{qqy})\omega^1-\tfrac{\gamma}{\alpha}\der\alpha \\
  &\Omega^2=-\tfrac{1}{6\alpha}F_{qq}(\tfrac{1}{2}\gamma^2+\tfrac{1}{3}F_q\gamma+K)\omega^4 \nonumber \\
   &+\tfrac{1}{6\alpha}(-\tfrac{1}{2}F_{qqq}\gamma^2+(\tfrac{1}{3}F_{qqq}F_{q}+F_{qqp})\gamma+F_{qqq}K-F_{qqy})\omega^2
   \nonumber \\ 
   &+\tfrac{1}{6\alpha}(-\tfrac{1}{2}F_{qqq}\gamma^3+(\tfrac{1}{6}F_{qq}^2+\tfrac{1}{3}F_{qqq}F_{q}+F_{qqp})\gamma^2
    +(F_{qq}K_q-F_{qqy}+F_{qqq}K)\gamma\nonumber \\
   &-\tfrac{1}{3}F_{qq}F_{qy}-F_{qq}K_p-\tfrac{1}{3}F_{qq}F_{q}K_{q}+\tfrac{1}{3}F_{qq}^2 K)\omega^1
    +\tfrac{1}{6\alpha}F_{qq}\der \gamma, \nonumber  
  \eeq
where $K$ denotes
\ben
 K=\tfrac{1}{6}(F_{qx}+pF_{qy}+qF_{qp}+FF_{qq})-\tfrac{1}{9}F_{q}^2-\tfrac{1}{2}F_p
\een
and $\omega^i$, $i=1,2,3,4$ are defined by the ODE via \eqref{omega}.
\end{theorem}
\noindent Exterior derivatives of the above invariant forms read
\beq\label{dtheta}
  &\der\theta^1 =\Omega^1\wedge\theta^1+\theta^4\wedge\theta^2,\nonumber \\
  &\der\theta^2 =\Omega^2\wedge\theta^1+a\theta^3\wedge\theta^2+b\theta^4\wedge\theta^2
    +\theta^4\wedge\theta^3,\nonumber \\
  &\der\theta^3 =\Omega^2\wedge\theta^2-\Omega^1\wedge\theta^3+(2-2c)\,\theta^3\wedge\theta^2
   +e\,\theta^4\wedge\theta^1+2b\,\theta^4\wedge\theta^3,\nonumber \\
  &\der\theta^4 =\Omega^1\wedge\theta^4+f\,\theta^4\wedge\theta^1+(c-2)\,\theta^4\wedge\theta^2
   +a\,\theta^4\wedge\theta^3,\\
  &\der\Omega^1 =(2c-2)\,\Omega^2\wedge\theta^1-\Omega^2\wedge\theta^4+g\,\theta^1\wedge\theta^2
   +h\,\theta^1\wedge\theta^3+k\,\theta^1\wedge\theta^4-f\,\theta^2\wedge\theta^4,\nonumber \\
  &\der\Omega^2 =\Omega^2\wedge\Omega^1-a\Omega^2\wedge\theta^3-b\,\Omega^2\wedge\theta^4+l\,\theta^1\wedge\theta^2
   +m\,\theta^1\wedge\theta^3+n\,\theta^1\wedge\theta^4 \nonumber\\
   &+r\,\theta^2\wedge\theta^3+s\,\theta^2\wedge\theta^4
   -f\,\theta^3\wedge\theta^4,\nonumber 
\eeq
where $a,b,c,e,f,g,h,k,l,m,n,r,s$ are functions on $\mathcal{P}$, which can be simply calculated due to formulae \eqref{theta}.
The simplest and the most symmetric case, when all the functions $a,b,c,e,f,g,h,k,l,m,n,r,s$ vanish, 
corresponds to the f.p. equivalence class of 
equation $$y'''=\tfrac{3}{2}\tfrac{{y''}^2}{y'}.$$ In this case, the 
manifold $\mathcal{P}$ is (locally) the Lie group
$SO(2,2)$ and the coframe
$(\theta^1,\theta^2,\theta^3,\theta^4,\Omega^1,\Omega^2)$ is a basis
of left invariant forms, which can be collected to the
$so(2,2)$-valued flat Cartan connection on $\mathcal{P}=SO(2,2)$. 
Since the Levi-Civita connection for the split signature metrics in
four dimensions also takes value in $so(2,2)$,
we ask under which conditions on f.p. equivalence classes of ODEs
(\ref{r1}), the equations \eqref{dtheta} may be interpreted as 
the structure equations for the Levi-Civita connection of a certain
4-dimensional split signature metric $G$.

\section{The construction of the metrics}
It is convenient to change the basis of 1-forms 
$\theta^1,\theta^2,\theta^3,\theta^4,\Omega^1,\Omega^2$ on $\mathcal{P}$ to
\begin{align}
 \tau^1&=2\theta^1+\theta^4, &  \tau^2&=\Omega^2, & \tau^3&=\Omega^2+2\theta^3, & \tau^4&=\theta^4, \label{tau}\\
 && \Gamma_1&=\Omega^1, & \Gamma_2&=\Omega^1+2\theta^2. && \nonumber 
\end{align}
After this change, equations (\ref{dtheta}) yield the formulae 
for the exterior differentials of
$\tau^1,\tau^2,\tau^3,\tau^4,\Gamma_1,\Gamma_2$. These are the 
formulae (\ref{dtau}) of the appendix. They can be
used to analyze the properties of the following bilinear tensor field
\be\label{metryka} 
 \wt{G}=\wt{G}_{ij}\tau^i\tau^j=2\tau^1\tau^2+2\tau^3\tau^4
\ee
on $\mathcal{P}$. The first question we ask here is the following:
under which conditions on $a,b,c,e,f,g,h,k,l,m,n,r,s$ the first four
of equations \eqref{dtau} may be identified with 
\ben \der\tau^i+\Gamma^i_{~j}\w\tau^j=0, \een
where the 1-forms $\Gamma^i_{~j}$, $ i,j=1,2,3,4 $ satisfy 
\ben
 \Gamma_{(ij)}=0, \qquad \textrm{and} \qquad   \Gamma_{ij}=\wt{G}_{ik}\Gamma^k_{~j}.
\een
This happens if and only if 
\be
c=0,\quad\quad l=0,\quad\quad r=0,\quad\quad s=0.\label{cond1}
\ee 
Now, we 
call 1-forms $\Gamma_1$, $\Gamma_2$ as {\it vertical} and 1-forms 
$\tau^1,\tau^2,\tau^3,\tau^4$ as {\it horizontal}.
To be able to interprete 
\ben R^i_{~j}=\der\Gamma^i_{~j}+\Gamma^i_{~k}\w\Gamma^k_{~j} \een
as a curvature, we have to require that 
it is horizontal, i.e contains no $\Gamma_1,\Gamma_2$ terms. This 
is equivalent to 
\be
m=0,\quad\quad a=0,\quad\quad g=0,\quad\quad f=-b.
\label{cond2}
\ee 
If these conditions are
satisfied then the exterior derivatives of \eqref{dtau} give also 
\be
b=0,\quad\quad h=0.\label{cond3}
\ee 
Concluding, having conditions (\ref{cond1}), (\ref{cond2}) and
(\ref{cond3}) satisfied, we have the following 
differentials of the coframe 
$(\theta^1,\theta^2,\theta^3\theta^4,\Gamma_1,\Gamma_2)$:
\begin{align}
 \der\tau^1=&\Gamma_1\w\tau^1,\nonumber  \\ 
 \der\tau^2=&-\Gamma_1\w\tau^2+\tfrac{1}{2}n\tau^1\w\tau^4, \nonumber \\
 \der\tau^3=&-\Gamma_2\w\tau^3+\left(\tfrac{1}{2}n-e\right)\tau^1\w\tau^4, \nonumber \\ 
 \der\tau^4=& \Gamma_2\w\tau^4, \label{dtauOK} \\ 
\der\Gamma_1=&\tau^1\w\tau^2+\tfrac{1}{2}k\tau^1\w\tau^4, \nonumber\\
\der\Gamma_2=&\tfrac{1}{2}k\tau^1\w\tau^4-\tau^3\w\tau^4, \nonumber 
\end{align}
and the following formulae for the matrix of 1-forms  
\ben
\Gamma^i_{~j}=\bma 
  -\Gamma_1 & 0 & 0 & 0 \\
  0 & \Gamma_1 & 0 & -\tfrac{1}{2}n\tau^1+(e-\tfrac{1}{2}n)\tau^4  \\                                       
   \tfrac{1}{2}n\tau^1-(e-\tfrac{1}{2}n)\tau^4 & 0 & \Gamma_2 & 0 \\
  0 & 0 & 0 & -\Gamma_2
\ema.\een
Moreover, introducing the frame of the vector fields $(X_1,X_2,X_3,X_4,
Y_1,Y_2)$ dual to the coframe $\tau^1,\ldots,\tau^4,\Gamma_1,\Gamma_2$
we get the following non-vanishing 2-forms $R^i_{~j}$:
\begin{align*}
 R^1_{~1}=&-\tau^1\w\tau^2-\tfrac{1}{2}k\tau^1\w\tau^4, \\
 R^2_{~2}=&\tau^1\w\tau^2+\tfrac{1}{2}k\tau^1\w\tau^4, \\
 R^2_{~4}=&\tfrac{1}{2}k\tau^1\w\tau^2+(\tfrac{1}{2}n_4+e_1-\tfrac{1}{2}n_1)\tau^1\w\tau^4-\tfrac{1}{2}k\tau^3\w\tau^4, \\  
 R^3_{~1}=&-\tfrac{1}{2}k\tau^1\w\tau^2-(\tfrac{1}{2}n_4+e_1-\tfrac{1}{2}n_1)\tau^1\w\tau^4+\tfrac{1}{2}k\tau^3\w\tau^4, \\ 
 R^3_{~3}=&\tfrac{1}{2}k\tau^1\w\tau^4-\tau^3\w\tau^4, \\
 R^4_{~4}=& -\tfrac{1}{2}k\tau^1\w\tau^4+\tau^3\w\tau^4.
\end{align*}
Here $f_i$ denotes $X_i(f)$. It further follows that 
$Ric_{ij}=R^k_{~ikj}$ satisfies  
\be Ric_{ij}=-\wt{G}_{ij}. \label{einstein}\ee  

These preparatory steps enable us to associate with each f.p. 
equivalence class of ODEs (\ref{r1}) satisfying conditions
(\ref{cond1})-(\ref{cond3}) a 4-manifold $\mathcal{M}$ 
equipped with a split signature Einstein metric $G$. 
This is done as follows.
\begin{itemize}
\item  The system (\ref{dtauOK})
  guarantees that the distribution $\mathcal{V}$ spanned by the vector fields
  $Y_1,Y_2$ is integrable. The leaf space of this foliation is
  4-dimensional and may be identified with $\mathcal{M}$. We also have
  the projection $\pi:\mathcal{P}\to\mathcal{M}$.
\item
The tensor field $\wt{G}$ is degenerate, $\wt{G}(Y_1,\cdot)=0$, 
$\wt{G}(Y_2,\cdot)=0$, along the leaves
of $\mathcal{V}$. Moreover, equations 
\eqref{dtauOK} imply that
\ben   
 L_{\scriptscriptstyle Y_1}\wt{G}=0, \qquad L_{\scriptscriptstyle Y_2}\wt{G}=0.
\een
Thus, $\wt{G}$ projects to a well defined split signature metric $G$
on $\mathcal{M}$.
\item
The Levi-Civita connection 1-form for $G$ and the curvature 2-form, 
pull-backed via $\pi^*$ to $\mathcal{P}$, 
identify with $\Gamma^i_{~j}$ and $R^i_{~j}$ respectively
\item 
Thus, due to equations (\ref{einstein}), the metric $G$ satisfies the Einstein
field equations with cosmological constant $\Lambda=-1$.
\end{itemize}
\label{pocz}

Below we find all functions $F=F(x,y,p,q)$ which 
solve conditions (\ref{cond1})-(\ref{cond3}). This will enable us 
to write down the explicit formulae for the Einstein 
metrics $G$ associated with 
the corresponding equations $y'''=F(x,y,y',y'')$.

The conditions $b=0$, $c=0$ in coordinates $x,y,p,q,\alpha,\gamma$ read 
\ben
F_{qp}+\tfrac{1}{3}F_{qq}+3K_{q}=0,\qquad 
F_{qqq}\gamma-F_{qqp}-\tfrac{1}{3}F_{qqq}F_{q}+\tfrac{1}{6}F_{qq}^2=0.
\een
The most general funtion $F(x,y,p,q)$ defining 3rd order ODEs satisfying these constraints is
\ben
F=\tfrac{3}{2}\frac{q^2}{p+\sigma(x,y)}+3\frac{\sigma_x(x,y)+p\sigma_y(x,y)}{p+\sigma(x,y)}q+\xi(x,y,p),
\een
where $\sigma,\xi$ are arbitrary functions of two and three varaibles,
respectively. Since the equations are considered modulo fiber preserving
transformations, we can put $\sigma=0$ by transformation $\bar{x}=x$
and $\bar{y}=\bar{y}(x,y)$ such that 
$\bar{y}_x=-\sigma(x,\bar{y}(x,y))$. Condition $l=0$ now becomes
\ben
p^3\xi_{ppp}-3p^2\xi_{pp}+6p\xi_p-6\xi=0,
\een
with the following general solution
\ben
 \xi=A(x,y)p^3+C(x,y)p^2+B(x,y)p.      
\een
Hence $F$ is given by 
\be\label{rdobre}
 F=\tfrac{3}{2}\frac{q^2}{p}+A(x,y)p^3+C(x,y)p^2+B(x,y)p.     
\ee
It further follows that it 
fulfills the remaining conditions $a=f=g=h=m=r=s=0$ and that 
\be\label{ccop}
k=-\frac{C}{4\alpha^2 p}\quad\quad n=\frac{C_y-zC-2A_x}{8\alpha^3
  p}\quad\quad e=\frac{1}{2}n+\frac{t C+2B_y-C_x}{16\alpha^3p^2} .
\ee

A straigthfowrward application of Theorem \ref{tw2} leads to the
following expressions for the `null coframe'
$(\tau^1,\tau^2,\tau^3,\tau^4)$:
\beqn
\tau^1&=&2\alpha~\der y\nonumber\\
\tau^2&=&(4\alpha)^{-1}~[~C~\der x+(2A -z^2)~\der y+2\der z~]\label{dcof}\\
\tau^3&=&(4\alpha p)^{-1}~[~-(t+2B)~\der x-C~\der y+2\der t ~]\nonumber\\
\tau^4&=&2\alpha p~\der x,\nonumber
\eeqn
where the new coordinates $z$ and $t$ are  
$$z=\frac{\gamma}{p}\quad\quad t=\frac{q}{p}+\gamma.$$
This brings $$\wt{G}=2(\tau^1\tau^2+\tau^3\tau^4)$$ on $\mathcal{P}$ to the form that depends only on
coordinates $(x,y,z,t)$. Thus, $\wt{G}$ projects to a well defined split signature 
metric 
\ben
G=-[t^2+2B(x,y)]\der x^2+2\der t\der x+[2A(x,y)-z^2]\der
 y^2+2\der z\der y 
\een
on a 4-manifold $\mathcal{M}$ parmetrized by
$(x,y,z,t)$. 

It follows from the construction that metric $G$ is
f.p. invariant. However, it does not yield all the f.p. information
about the corresponding ODE. It is clear, since the function $C$ which is
proportional to the f. p. Cartan's invariant $k$ of (\ref{dtauOK}), 
is not appearing in the
metric $G$. From the point of view of the metric, function $C$ represents a
`null rotation' of coframe $(\tau^i)$. Thus it is not a geometric
quantity. Therefore $G$, although f.p. invariant, can not distinguish 
between various f.p. nonequivalent classes of equations 
such as, for example, those with
$C\equiv 0$ and $C\neq 0$. To fully distinguish all non-equivalent ODEs
with (\ref{rdobre}) one needs additional structure than the metric
$G$. This structure is only fully described by the bundle
$\pi:\mathcal{P}\to\mathcal{M}$ together with the coframe 
$(\tau^1,\tau^2,\tau^3,\tau^4,\Gamma_1,\Gamma_2)$ of
(\ref{dtauOK}) on $\mathcal{P}$. An alternative description, 
more in the spirit of the split signature metric $G$, is 
presented in section 5.

Now, equations (\ref{einstein}) imply that the metric $G$ is 
Einstein with cosmological consatnt $\Lambda=-1$. The anti-selfdual part of
its Weyl tensor is always of Petrov-Penrose 
type D. The selfdual Weyl tensor is of type II, if the functions $A$ and $B$
are generic. If $A=A(y)$ and  
$B=B(x)$ the selfdual Weyl tensor degenerates to a tensor of type D. 
Summing up we have following theorem.
\begin{theorem}\label{tw1}
Third order ODE 
\ben
 y'''= \tfrac{3}{2}\frac{{y''}^2}{y'}+A(x,y){y'}^3+C(x,y){y'}^2+B(x,y)y'  
\een
defines, by virtue of Cartan's equivalence method, a 4-dimensional split signature metric
\ben
G=-[t^2+2B(x,y)]\der x^2+2\der t\der x+[2A(x,y)-z^2]\der
 y^2+2\der z\der y 
\een
which is Einstein
\ben
Ric(G)=-G
\een
and has Weyl tensor $W=W^{\scriptscriptstyle ASD}+W^{\scriptscriptstyle SD}$ of Petrov type D+II, with the exception of the case 
$A=A(y)$, $B=B(x)$, when it is of type D+D. The  metric $G$ is
invariant with respect to f.p. transformations of the variables of the ODE. 
\end{theorem}

\section{Uniqueness of the metrics}  
In this section we prove the following theorem.
\begin{theorem}\label{tw3}
The metrics of theorem \ref{tw1} are the unique family of metrics $G$,
which are defined by f.p. equivalence classes of 3rd order ODEs and
satisfy the following three conditions.
\begin{itemize}
 \item The metrics are split signature, Einstein: $Ric(G)=-G$, and each of them is defined on 4-d manifold 
   $\mathcal{M}$, which is the base of the fibration
   $\pi:\mathcal{P}\to\mathcal{M}$.
\item The family contains a metric corresponding to equation $y'''=\tfrac{3}{2}\tfrac{{y''}^2}{y'}$.
 \item The tensor 
\ben
 \wt{G}=\pi^*G=\mu_{ij}\theta^i\theta^j+\nu_{iA}\theta^i\Omega^A+\rho_{AB}\Omega^A\Omega^B,
 \een
on $\mathcal{P}$, when expressed by 
the invariant coframe $(\theta^i,\Omega^A)$ 
associated with the respective f.p. equivalence class,  
 has the coefficients $\mu_{ij}$, $\nu_{iA}$, $\rho_{AB}$; $i,j=1,\ldots,4$; $A,B=1,2$ \emph{constant and the same} for all classes of the ODEs 
 for which  $G$ is defined.  
\end{itemize}
\end{theorem}
\noindent
To prove the theorem, it is enough to show the 
uniqueness of $G$ in the simplest case of equation $y'''=\tfrac{3}{2}\tfrac{{y''}^2}{y'}$,
and to repeat the calculations from pp. \pageref{pocz} --
\pageref{rdobre} 
of Section \ref{pocz} for a generic equation.
The following trivial proposition holds.
\begin{proposition}\label{prop}
Let $\wt{G}$ be a bilinear symmetric form of signature $(++--00)$ on 
$\mathcal{P}$, such that for a vector field $N$ 
\be\label{killing}
if\qquad \wt{G}(N,\cdot)=0\qquad then\qquad L_{\scriptscriptstyle N}\wt{G}=0. 
\ee
A distribution spanned by such vector fields $N$ 
is integrable and defines a 
4-d manifold $\mathcal{M}$ as a space of its integral
leaves. There exists exactly one bilinear form $G$ on $\mathcal{M}$
with the property $\pi^{\ast}G=\wt{G}$, where
$\pi:\mathcal{P}\to\mathcal{M}$ is the canonical projection
assigning a point of $\mathcal{M}$ to an integral leave of the
distribution.
\end{proposition}
\noindent
Our aim now is to find all the metrics $\wt{G}$ of proposition \ref{prop} which, when expressed by the coframe $\theta^i,\Omega^A$ 
(or, equivalently, by $\tau^i,\Gamma_A$), have constant coefficients. Let us consider the simplest case, corresponding to equation 
$y'''=\tfrac{3}{2}\tfrac{{y''}^2}{y'}$, for which all the invariant
functions appearing in \eqref{dtheta} and \eqref{dtau}
vanish. $\mathcal{P}$ is now the Lie group 
$SO(2,2)$, $\wt{G}$ is a form on Lie algebra $so(2,2)$, the
distribution spanned by the degenerate fields $N$ is a 2-d 
subalgebra $\gt{h}\subset so(2,2)$. 
Finding $\wt{G}$ is now a purely algebraic problem. 
In our case the basis $(\tau^i,\Gamma_A)$ satisfies
\begin{align}
 \der\tau^1&=\Gamma_1\w\tau^1, &  \der\tau^3&=-\Gamma_2\w\tau^3, \nonumber \\
 \der\tau^2&=-\Gamma_1\w\tau^2, &  \der\tau^4&=\Gamma_2\w\tau^4, \label{dtau_alg} \\
 \der\Gamma_1&=\tau^1\w\tau^2,  &  \der\Gamma_2&=\tau^4\w\tau^3, \nonumber
 \end{align}
which agrees with a decomposition $so(2,2)=so(1,2)\oplus so(1,2)$. 
A group of transformations preserving equations \eqref{dtau_alg} is $O(1,2)\times O(1,2)$, 
that is the intersection of the orthogonal group $O(2,4)$ preserving
the Killing form $\kappa$ of $so(2,2)$ and the group 
$GL(3)\times GL(3)$ preserving the decomposition $so(2,2)=so(1,2)\oplus so(1,2)$. 
Each coframe $(\wt{\tau}^i$, $\wt{\Gamma}_A)$, satisfying \eqref{dtau_alg} is obtained by a linear transformation:
\be\label{change_b}
 \bma \wt{\tau}^1 \\  \wt{\tau}^2 \\ \wt{\Gamma}_1 \ema= A \bma \tau^1 \\  \tau^2 \\ \Gamma_1 \ema, \qquad
  \bma \wt{\tau}^3 \\  \wt{\tau}^4 \\ \wt{\Gamma}_2 \ema= B \bma \tau^3 \\  \tau^4 \\ \Gamma_2 \ema 
  \qquad A,B\in O(1,2). 
\ee  
We use transformations  \eqref{change_b} to obtain the most
convenient form of the basis
$(N_1,N_2)$ of the subalgebra $\gt{h}\subset so(2,2)$. 
We write down the metric $\wt{G}$ in the corresponding 
coframe $(\wt{\tau}^1,\wt{\tau}^2,\wt{\tau}^3,\wt{\tau}^4,\wt{\Gamma}_1,\wt{\Gamma}_2)$
 and impose conditions
(\ref{killing}). 
This conditions imply that the most general form of the metric is 
$\wt{G}=2u\wt{\tau}^1\wt{\tau}^2+2v\wt{\tau}^3\wt{\tau}^4,$ where $u,v$
are two real parameters. In such case, $[N_1,N_2]=0$ and 
$\kappa(N_1,N_1)<0$, $\kappa(N_2,N_2)<0$. 
When written in terms of the coframe  $(\tau^i,\Gamma_A)$, $\wt{G}$ involves six real parameters 
$u,v,\mu,\phi,\nu,\psi$, however it appears, that only parameters $u$ and $v$ are essential; 
different choices of $\mu,\phi,\nu,\psi$ define different degenerate distributions spanned by
$N_1,N_2$ and hence spaces $\mathcal{M}$ are different, but metrics $G$ on them are isometric. Thus we can choose 
$\wt{G}=2u\tau^1\tau^2+2v\tau^3\tau^4$.
Computing $\wt{G}$ for $F=\tfrac{3}{2}\frac{q^2}{p}$, we have, in a
suitable coordinate system $(x,y,z,t)$, 
\ben 
 G=-v[t^2+2B(x,y)]\der x^2+2v\der t\der x+u[2A(x,y)-z^2]\der
 y^2+2u\der z\der y.  
\een 
Parameters $u,v$ can be also fixed, if we demand $G$ to be Einstein with cosmological constant $\Lambda=-1$. 
This is only possible if $u=1$, $v=1$. The the tensor field  $\wt{G}$
defined in this way is unique and has the form 
\ben
 \wt{G}=2\tau^1\tau^2+2\tau^3\tau^4=2\Omega^2(2\theta^1+\theta^4)+2\theta^4(2\theta^3+\Omega^2).
\een
This formula is used in the generic case explaining our choice of the
coframe (\ref{tau}) and the metric (\ref{metryka}). This finishes the proof of theorem \ref{tw3}.

\section{The Cartan connection and the distinguished class of ODEs}
Here we provide an alternative description of the f.p. equivalence class of
third order ODEs corresponding to $F=F(x,y,p,q)$ of (\ref{rdobre}). We
consider a 4-dimensional manifold $\mathcal{M}$ 
parametrized by $(x,y,z,t)$. Then the geometry of a f.p. equivalence
class of ODEs (\ref{rdobre}) is in one to one correspondence with the geometry
of a class of coframes
\beq
\tau^1_0&=&\der y\nonumber\\
\tau^2_0&=&\frac{1}{2}~[~C~\der x+(2A -z^2)~\der y+2\der z~]\nonumber\\
\tau^3_0&=&\frac{1}{2}~[~-(t+2B)~\der x-C~\der y+2\der t ~]\label{traf1}\\
\tau^4_0&=&\der x,\nonumber
\eeq
on $\mathcal{M}$ given modulo
a special $SO(2,2)$ transformation 
\be
\tau^i_0\mapsto\tau^i= h^i_{~j}\tau^j_0,\quad\quad {\rm where}\quad\quad 
(h^i_{~j})=
\bma 2\alpha&0&0&0\\
0&(2\alpha)^{-1}&0&0\\
0&0&(2\alpha p)^{-1}&0\\
0&0&0&2\alpha p
\ema.\label{traf}
\ee

The Cartan 
equivalence method applied to the question if two coframes
(\ref{traf1}) are transformable to each other via (\ref{traf}) 
gives the full system of invariants of this geometry. These invariants
consist of (i) a fibration $\pi:\mathcal{P}\to\mathcal{M}$ of Section
3, which now becomes a Cartan bundle
$\mathcal{H}\to\mathcal{P}\to\mathcal{M}$ with the 2-dimensional
structure group $\mathcal{H}$ generated by $h^i_{~j}$, and (ii)  of an $so(2,2)$-valued Cartan connection $\omega$ described by
the coframe $(\tau^1,\tau^2,\tau^3,\tau^4,\Gamma_1,\Gamma_2)$ of
(\ref{dtauOK}) on $\mathcal{P}$. Explicitely, the connection $\omega$ is
given by
$$
\omega^i_{~j}=\bma- \frac{1}{2}(\Gamma_1+\Gamma_2+\tau^4)&0&\tau^1&-\frac{1}{2}\tau^4\\
0&\frac{1}{2}(\Gamma_1+\Gamma_2+\tau^4)&-\Gamma_2+\tau^3-\frac{1}{2}\tau^4&-\frac{1}{2}\tau^2\\
\frac{1}{2}\tau^2&\frac{1}{2}\tau^4&\frac{1}{2}(\Gamma_1-\Gamma_2-\tau^4)&0\\
\Gamma_2-\tau^3+\frac{1}{2}\tau^4&-\tau^1&0&\frac{1}{2}(-\Gamma_1+\Gamma_2+\tau^4)
\ema.
$$
To see that this is an $so(2,2)$ connection it is enough to note that
 $g_{ij}\omega^j_{~k}+g_{kj}\omega^k_{~i}=0$ with the matrix $g_{ij}$
given by 
$$
g_{ij}=\bma 0&1&0&0\\
1&0&0&0\\
0&0&0&1\\
0&0&1&0
\ema.
$$
Now, equations (\ref{dtauOK}) are interpreted as the requirement that
the curvature $$\Omega=\der\omega+\omega\wedge\omega$$ of this connection
$\omega$ has a very simple form
$$\Omega=\bma -\frac{1}{2}k&0&0&0\\
0&\frac{1}{2}k& \frac{1}{2}(-k+n-2e)&-\frac{1}{4}n\\
\frac{1}{4}n&0&0&0\\
\frac{1}{2}(k-n+2e)&0&0&0
\ema\tau^1\wedge\tau^4,
$$
where $n,e$ and $k$ are given by (\ref{ccop}). The connection $\omega$
and its curvature $\Omega$ yields all the f.p. information of the
equation corresponding to (\ref{rdobre}). In particular, all the equations
with $k=n=e=0$ are f.p. equivalent, all having the vanishing curvature of
their Cartan connection $\omega$.

It is interesting to search for a split signature 4-metric
$H$ for which the connection $\omega$ is the 
Levi-Civita connection. The general form of such metric is 
$$H=g_{ij}T^iT^j,$$ where $(T^1,T^2,T^3,T^4)$ are four linearly
independent 1-forms on $\mathcal{P}$ which staisfy  
\be
{\rm d}T^i+\omega^i_{~j}\wedge T^j=0.\label{bianchi}
\ee
Thus, for such $H$ to exist, the 1-forms $(T^1,T^2,T^3,T^4)$ must also 
satisfy the integrability conditions of (\ref{bianchi}),  
$$\Omega^i_{~j}\wedge T^j=0,$$
which are just the Bianchi identities for $\omega$ to be the
Levi-Civita connection of metric $H$. 
These identities provide severe algebraic constraints on the possible
solutions $(T^i)$. Using them, under the assumption that
$C(x,y)\neq 0$ in the considered region of $\mathcal{P}$,  
we found all $(T^i)$s satisfying (\ref{bianchi}). Thus, with every
triple $C\neq 0,A,B$  corresponding to an ODE given by $F$ of
(\ref{rdobre}), we were able to find a split signature metric $H$ for
which connection $\omega$ is the Levi-Civita connection. Surprisingly,
given $A, B$ and $C\neq 0$ the general solution for $(T^i)$ involves 
four {\it free} real functions. Two of these functions 
depend on 6 variables and the other two depend on 2 variables. Thus, each
f.p. equivalence class of ODEs representd by $F$ of (\ref{rdobre})
defines a large family of split signature metrics $H$ for which
$\omega$ is the Levi-Civita connection\footnote{The 4-manifold on which each
of these metrics resides is the leaf space of the 2-dimensional
integrable distribution on $\mathcal{P}$ which anihilates forms $(T^1,T^2,T^3,T^4)$.}. 
Writing down the explicit formulae for these metrics is easy, but we
do not present them here, due to their ugliness and due to the fact
that, regardless of the choice of the four free functions, they never 
satisfy the Einstein equations. The proof of this last fact is based
on lengthy calculations using the explicit forms of the general
solutions for $(T^i)$. 

\section*{Appendix}
In this appendix we give the formulae for the differentials of the
transformed Cartan invariant coframe
$(\tau^1,\tau^2,\tau^3,\tau^4,\Gamma_1,\Gamma_2)$ on
$\mathcal{P}$. These are:
\begin{subequations}\label{dtau}
\begin{align}
  \der\tau^1=&\Gamma_1\w\tau^1+\tfrac{1}{2}c\Gamma_1\w\tau^4-\tfrac{1}{2}c\Gamma_2\w\tau^4
   +\tfrac{1}{2}\,f\tau^4\w\tau^1-\tfrac{1}{2}a\tau^4\w\tau^2 \label{dtau1} \\ 
   &+\tfrac{1}{2}a\tau^4\w\tau^3,\nonumber \\
 \der\tau^2=& \tfrac{1}{4}l\Gamma_1\w\tau^1+(\tfrac{1}{4}r-1)\Gamma_1\w\tau^2  -\tfrac{1}{4}r\Gamma_1\w\tau^3
  -\left(\tfrac{1}{4}l+\tfrac{1}{2}s\right)\Gamma_1\w\tau^4  \label{dtau2} \\
  &-\tfrac{1}{4}l\Gamma_2\w\tau^1 -\tfrac{1}{4}r\Gamma_2\w\tau^2  +\tfrac{1}{4}r\Gamma_2\w\tau^3 
  +\left(\tfrac{1}{4}l+\tfrac{1}{2}s\right)\Gamma_2\w\tau^4  \nonumber\\
   & +\tfrac{1}{4}m\tau^2\w\tau^1-\tfrac{1}{4}m\tau^3\w\tau^1 -\tfrac{1}{2}n\tau^4\w\tau^1  +\tfrac{1}{2}a\tau^3\w\tau^2 \nonumber \\
  &+(\tfrac{1}{4}m-\tfrac{1}{2}f+b)\tau^4\w\tau^2 +\left(\tfrac{1}{2}f-\tfrac{1}{4}m\right)\tau^4\w\tau^3,  \nonumber \\
 \der\tau^3=&\tfrac{1}{4}l\Gamma_1\w\tau^1  +\left(c+\tfrac{1}{4}r\right)\Gamma_1\w\tau^2  -\left(c+\tfrac{1}{4}r \right)\Gamma_1\w\tau^3  
  -\left(\tfrac{1}{4}l+\tfrac{1}{2}s \right)\Gamma_1\w\tau^4  \label{dtau3}\\ 
  &+\tfrac{1}{4}l\Gamma_2\w\tau^1  -\left(c+\tfrac{1}{4}r\right)\Gamma_2\w\tau^2  +\left(c+\tfrac{1}{4}r-1\right)\Gamma_2\w\tau^3 \nonumber\\
  &+\left(\tfrac{1}{4}l+\tfrac{1}{2}s \right)\Gamma_{2}\w\tau^{4} +\tfrac{1}{4}m\tau^2\w\tau^1  
  -\tfrac{1}{4}m\tau^3\w\tau^1  +\left(e-\tfrac{1}{2}n \right) \tau^4\w\tau^1 \nonumber\\ 
  &+\tfrac{1}{2}a\tau^3\w\tau^2  
  +\left(\tfrac{1}{4}m-b-\tfrac{1}{2}f\right)\tau^4\w\tau^2 +\left(2b+\tfrac{1}{2}f-\tfrac{1}{4}m\right)\tau^4\w\tau^3, \nonumber\\ 
 \der\tau^4=& +\tfrac{1}{2}c\Gamma_1\w\tau^4 +\left(1-\tfrac{1}{2}c\right)\Gamma_2\w\tau^4
    +\tfrac{1}{2}f\tau^4\w\tau^1 -\tfrac{1}{2}a\tau^4\w\tau^2  \label{dtau4} \\
    &+\tfrac{1}{2}a\tau^4\w\tau^3, \nonumber \\
\der\Gamma_1=& \tfrac{1}{4}g\Gamma_1\w\tau^1+\left(\tfrac{1}{2}f-\tfrac{1}{4}g\right)\Gamma_1\w\tau^4 
   -\tfrac{1}{4}g\Gamma_2\w\tau^1 +\left(\tfrac{1}{4}g-\tfrac{1}{2}f\right)\Gamma_2\w\tau^4 \label{dGamma1} \\
  &+\left( \tfrac{1}{4}h+c-1 \right)\tau^2\w\tau^1+-\tfrac{1}{4}h\tau^3\w\tau^1-\tfrac{1}{2}k\tau^4\w\tau^1  \nonumber \\&
  +\left( \tfrac{1}{4}h+c\right)\tau^4\w\tau^2  -\tfrac{1}{4}h\tau^4\w\tau^3, \nonumber\\
\der\Gamma_2=&\tfrac{1}{4}g\Gamma_1\w\tau^1-\tfrac{1}{2}a\Gamma_1\w\tau^2 +\tfrac{1}{2}a \Gamma_1\w\tau^3  
  +\left(b+\tfrac{1}{2}f-\tfrac{1}{4}g\right)\Gamma_1\w\tau^4 \label{dGamma2} \\ 
  &-\tfrac{1}{4}g\Gamma_2\w\tau^1  +\tfrac{1}{2}a\Gamma_2\w\tau^2  
  -\tfrac{1}{2}a\Gamma_2\w\tau^3 +\left(\tfrac{1}{4}g-b-\tfrac{1}{2}f\right)\Gamma_2\w\tau^4  \nonumber\\
  & +\left( \tfrac{1}{4}h+c \right)\tau^2\w\tau^1  -\tfrac{1}{4}h\tau^3\w\tau^1  -\tfrac{1}{2}k\tau^4\w\tau^1 
  +\left(\tfrac{1}{4}h+c \right)\tau^4\w\tau^2 \nonumber\\ 
  &+\left( 1-\tfrac{1}{4}h \right)\tau^4\w\tau^3. \nonumber 
\end{align}
\end{subequations}

\end{document}